\documentclass[reqno]{amsart}
\usepackage{amsfonts,graphicx,float}

\restylefloat{figure}

\makeatletter                                                  %
\def\th@plain{\slshape}                                        %
\makeatother                                                   %

\newcommand{\oi}{[0,1]}

\newcommand{\Zbb}{\mathbb{Z}}
\newcommand{\Qbb}{\mathbb{Q}}
\newcommand{\Rbb}{\mathbb{R}}
\newcommand{\Cbb}{\mathbb{C}}

\newcommand{\llgroup}{$\ell$-group}
\newcommand{\llgroups}{$\ell$-groups}

\newcommand{\llideals}{$\ell$-ideals}

\newcommand{\pfrak}{\mathfrak{p}}

\newcommand{\rfrak}{\mathfrak{r}}

\newcommand{\Luk}{\L ukasiewicz}

\newcommand{\newword}[1]{\textsl{#1}}

\newcommand{\vect}[3]{#1_#2,\ldots ,#1_#3}

\newcommand{\abs}[1]{\lvert#1\rvert}
\newcommand{\uniform}[1]{\lVert#1\rVert_\infty}
\newcommand{\angles}[1]{\langle #1 \rangle}

\DeclareMathSymbol{\upharpoonright}{\mathrel}{AMSa}{"16}
\let\restriction\upharpoonright
\DeclareMathSymbol{\nmid}{\mathrel}{AMSb}{"2D}

\DeclareMathOperator{\MaxSpec}{MaxSpec}
\DeclareMathOperator{\den}{den}

\DeclareMathOperator{\Free}{Free}
\DeclareMathOperator{\mmod}{mod}
\DeclareMathOperator{\supp}{supp}

\DeclareMathOperator{\SL}{SL}

\DeclareMathOperator{\Mat}{Mat}

\theoremstyle{plain}
\newtheorem{theorem}{Theorem}[section]
\newtheorem{lemma}[theorem]{Lemma}
\newtheorem{proposition}[theorem]{Proposition}

\theoremstyle{definition}
\newtheorem{definition}[theorem]{Definition}
\newtheorem{example}[theorem]{Example}

\begin{document}

\bibliographystyle{plain}

\sloppy

\title[Invariant measures]{Invariant measures\\
in free MV-algebras}

\author[G. Panti]{Giovanni Panti}
\address{Department of Mathematics\\
University of Udine\\
via delle Scienze 208\\
33100 Udine, Italy}
\email{panti@dimi.uniud.it}

\begin{abstract}
MV-algebras can be viewed either as the Lindenbaum algebras of \Luk\ infinite-valued logic, or as unit intervals of lattice-ordered abelian groups in which a strong order unit has been fixed.
The free $n$-generated MV-algebra $\Free_n$ is representable as an algebra of continuous piecewise-linear functions with integer coefficients over the unit cube $\oi^n$. The maximal spectrum of
$\Free_n$ is canonically homeomorphic to $\oi^n$, and the automorphisms of the algebra are in $1$--$1$ correspondence with the pwl homeomorphisms with integer coefficients of the unit cube. In this paper we prove that the only probability measure on $\oi^n$ which is null on underdimensioned $0$-sets and is invariant under the group of all such homeomorphisms is the Lebesgue measure. From the viewpoint of lattice-ordered abelian groups, this fact means that, in relevant cases, fixing an automorphism-invariant strong unit implies fixing a distinguished probability measure on the maximal spectrum. From the viewpoint of algebraic logic, it means that the only automorphism-invariant truth averaging process that detects pseudotrue propositions is the integral with respect to Lebesgue measure.
\end{abstract}

\keywords{MV-algebras, state, automorphism-invariance, piecewise-linear homeomorphisms}

\thanks{\emph{2000 Math.~Subj.~Class.}: 06D35; 37A05}

\maketitle

\section{Preliminaries}

An \newword{MV-algebra} is an algebra $(A,\oplus,\neg,0)$ such that $(A,\oplus,0)$ is a commutative monoid and the identities $\neg\neg f=f$, $f\oplus\neg 0=\neg 0$, and $\neg(\neg f\oplus g)\oplus g=\neg(\neg g\oplus f)\oplus f$ are satisfied. MV-algebras can be viewed either as the Lindenbaum algebras of \Luk\ infinite-valued logic, or as unit intervals of lattice-ordered abelian groups (\llgroups) in which a strong order unit has been fixed.
See~\cite{bkw}, \cite{beynon77}, \cite{glass99} for \llgroups, and \cite{mundicijfa}, \cite{mundici95}, \cite{CignoliOttavianoMundici00} for MV-algebras and \Luk\ logic.

We recall that a \newword{strong unit} in an \llgroup\ $G$ is a positive element $u$ of $G$ such that for every $g\in G$ there exists a positive integer $n$ for which $g\le nu$. The unit interval $\Gamma(G,u)=\{g\in G:0\le g\le u\}$ is then an MV-algebra under the operations $f\oplus g=(f+g)\land u$, $\neg f=u-f$, $0=0_G$, and the functor $\Gamma$ is an equivalence between the category of \llgroups\ with a distinguished strong unit and the category of MV-algebras~\cite{mundicijfa}. A \newword{state} on $(G,u)$ is a group homomorphism $m:G\to\Rbb$ which is \newword{positive} ($f\ge 0$ implies $m(f)\ge 0$) and \newword{normalized}
($m(u)=1$)~\cite[Chapter~4]{goodearl86}. The same definition, recast in terms of MV-algebras, amounts to the following: a \newword{state} on the MV-algebra $A$ is a function $m:A\to\oi$ such that $m(0)=0$, $m(\neg0)=1$, and $m(f\oplus g)=m(f)+m(g)$ provided that $\neg(\neg f\oplus\neg g)=0$~\cite{mundici95}. If $A$ is viewed as the Lindenbaum algebra of some theory in \Luk\ logic, then $m$ is a function assigning an ``average truth-value'' to the elements of $A$, i.e., to the propositional formulas modulo the theory.

The $\Gamma$ functor induces a canonical bijection between the states of $(G,u)$ and those of $A=\Gamma(G,u)$ and a homeomorphism between 
the maximal spectrum of $(G,u)$ and $\MaxSpec A$. 
We recall that the maximal spectrum of $(G,u)$ is the set of maximal \llideals\ of $G$, while $\MaxSpec A$ is the set of maximal ideals of $A$ (i.e., kernels of homomorphisms from $A$ to $\Gamma(\Rbb,1)$);
both sets are equipped with the Zariski topology~\cite[Chapitre~10]{bkw}.
We will formulate our results mainly in terms of MV-algebras, leaving to the reader their straightforward translation in the language of \llgroups\ with strong unit.

The set of all states of $A$ is a compact convex subset of $\oi^A$, where the latter is given the product topology, and the subspace $X$ of extremal states (a state $m$ is \newword{extremal} if it cannot be expressed as $m=tp+(1-t)q$, with $p,q$ distinct states and $0<t<1$)
is homeomorphic to $\MaxSpec A$. This homeomorphism is canonical: if $\pfrak$ is a maximal ideal, then there is a unique embedding of MV-algebras $A/\pfrak\to\Gamma(\Rbb,1)=\oi$, and the extremal state $p$ corresponding to $\pfrak$ is the composition of this embedding with the quotient map $A\to A/\pfrak$. We identify $X$ with $\MaxSpec A$, and we write $f(p)$ for the image of $f/\pfrak$ in $\oi$. As a function, $f$ is an element of $C(X)$ (the Banach vector lattice of real-valued continuous functions on $X$ endowed with the uniform norm $\uniform{\phantom{f}}$), and the resulting map $A\to\Gamma(C(X),1)$ (here $1$ is the characteristic function of~$X$)
is a homomorphism of MV-algebras whose kernel is the \newword{radical ideal} $\rfrak=\bigcap\MaxSpec A$ of $A$.

\begin{proposition}
Let\label{ref1} $A$ be an MV-algebra, $X=\MaxSpec A$, $P(X)$ the set of all regular Borel probability measures on $X$. Then the states on $A$ are in 1-1 correspondence with the elements of $P(X)$.
\end{proposition}
\begin{proof}
The states of $A$ are in 1-1 correspondence with the states of $A/\rfrak$~\cite[Proposition~3.1(3)]{mundici95}. We may then assume that $A$ has $0$ radical, and identify $A$ with a separating subalgebra of $\Gamma(C(X),1)$. Obviously every Borel probability measure $\mu$ on $X$ induces a state $m$ on $A$ via Lebesgue integration:
$$
m(f)=\int_X f\,d\mu.
$$
We must show that the map $\mu\mapsto m$, when restricted to the set of regular measures, is a bijection. Let then $m$ be a state on $A$, which uniquely corresponds to a state ---again denoted by $m$--- on the \llgroup\ with strong unit $(G,1)\subseteq(C(X),1)$
enveloping~$A$ \cite[Theorem~2.4]{mundici95}.
Since $G$ is torsion-free, it embeds in its divisible hull $H=\{qf:q\in\Qbb, f\in G\}$, which is a sub-vector lattice of $C(X)$ over the rationals. One sees easily that $m$ extends uniquely to a state on $(H,1)$, via $m(qf)=q\,m(f)$. By the lattice version of the Stone-Weierstrass Theorem, $H$ is dense in $C(X)$. Since $m$ is positive, it is bounded ($\uniform{f}\le 1$ $\implies$ $-1\le f\le 1$ $\implies$ $m(-1)\le m(f)\le m(1)$ $\implies$ $\abs{m(f)}\le m(1)$), whence continuous. Therefore $m$ can be extended uniquely to a state on $(C(X),1)$ via $m(f)=\lim m(h_n)$, whenever $f\in C(X)$ is the uniform limit of elements $h_n\in H$. Note that $m$ is still positive on $C(X)$, since every $0\le f\in C(X)$ can be written as the limit of positive elements of $H$. The proof is now completed by an application of the Riesz Representation Theorem~\cite[Theorem~2.14]{rudin87}, which guarantees the existence of a unique $\mu\in P(X)$ that induces~$m$.
\end{proof}

Let us note that if every open set in $X$ is $\sigma$-compact (i.e., a countable union of compact sets), then every Borel measure on $X$ is regular~\cite[Theorem~2,18]{rudin87}. This applies, in particular, when $A$ is a countable MV-algebra. Indeed, in this case $X$ is second countable, whence is metrizable and has $\sigma$-compact open sets. For countable $A$, we can then remove the word ``regular'' in the statement of Proposition~\ref{ref1}.

If $A,B$ are MV-algebras (not necessarily having $0$ radical), $X=\MaxSpec A$, $Y=\MaxSpec B$, $\sigma:A\to B$ a homomorphism, then the \newword{dual} $S:Y\to X$ of $\sigma$ is defined by $S(m)=m\circ\sigma$ (at the level of extremal states), or by $S(\pfrak)=\sigma^{-1}[\pfrak]$ (at the level of maximal ideals).
Recall that if $S:Y\to X$ is a Borel measurable map between topological spaces and $\mu$ is a Borel measure on $Y$, then the \newword{push-forward} of $\mu$ by $S$ is the Borel measure $S_*\mu$ on $X$ defined by $(S_*\mu)(A)=\mu(S^{-1}A)$.
If now $\mu\in P(Y)$ corresponds to $m$ as in Proposition~\ref{ref1}, then $S(m)$ corresponds to $S_*\mu$. Indeed, for every $f\in A$, we have $\sigma(f)=f\circ S$ as a function on $Y$, and therefore
$$
[S(m)](f)=m(\sigma(f))=\int_Y\sigma(f)\,d\mu=\int_Yf\circ S\,d\mu=
\int_Xf\,d(S_*\mu).
$$
If in particular $\sigma$ is an endomorphism of $A=B$, then $m$ is $\sigma$-invariant iff $(X,\mu,S)$ is a measure-theoretic dynamical system.

\section{The main result}

Fix $n\ge 1$, and let $\Free_n$ be the free MV-algebra over $n$ generators. It is well known that $\Free_n$ has $0$ radical and that $X=\MaxSpec\Free_n$ is homeomorphic to the $n$-cube $\oi^n$; see~\cite{mundicijfa}, \cite{mundicimn}, \cite{pantilu} for all unproved claims on free MV-algebras.
A \newword{rational polyhedral complex} is a finite set $W$ of compact convex polyhedrons such that: (1) every polyhedron is contained in $\oi^n$, and all its vertices have rational coordinates; (2) if $C\in W$ and $D$ is a face of $C$, then $D\in W$; (3) every two polyhedra intersect in a common face.
We write $\abs{W}$ for the union of all elements of $W$. The range of the embedding $\Free_n\to\Gamma(C(X),1)$ described before Proposition~\ref{ref1}
is then the set of all \newword{McNaughton functions on $\oi^n$}, i.e., those continuous functions $f:\oi^n\to\oi$ for which there exists a complex $W$
as above and affine linear functions $F_j(\bar x)=a_{j1}x_1+\cdots+a_{jn}x_n+a_{j(n+1)}\in\Zbb[\vect x1n]$, in 1-1 correspondence with the $n$-dimensional polyhedra $C_j$ of the complex, such that $\abs{W}=\oi^n$ and $f\restriction C_j=F_j$ for each $j$.
The $i$th projection $x_i:\oi^n\to\oi$ is the $i$th free generator of $\Free_n$. 
If $\sigma$ is an automorphism of $\Free_n$
then, upon the identification of $\MaxSpec\Free_n$ with $\oi^n$, its dual is the homeomorphism $S:\oi^n\to\oi^n$ defined by $S(p)_i=(\sigma(x_i))(p)$. 
This is equivalent to saying that there exists a complex $W$ as above and matrices
$P_j\in\Mat_{n+1}(\Zbb)$, in 1-1 correspondence with the $n$-dimensional polyhedra $C_j$ of $W$, such that $\abs{W}=\oi^n$ and:
\begin{enumerate}
\item every $P_j$ has last row $(0\ldots0\,1)$;
\item $P_j$ expresses $S\restriction C_j$ in homogeneous coordinates (i.e., if $p=(\vect r1n)\in C_j$ and
$S(p)=(\vect s1n)$, then $(s_1\ldots s_n\,1)^{tr}=
P_j(r_1\ldots r_n\,1)^{tr}$);
\item all $P_j$ have determinant $+1$, or they all have determinant $-1$.
\end{enumerate}
We call these duals \newword{McNaughton homeomorphisms} of the unit cube. Note that each of $\sigma$ and $S$ determines the other; indeed, given $S$ as above, $\sigma$ is the automorphism that maps $x_i$ to $x_i\circ S$.

Throughout this paper $\lambda$ (or $\lambda^n$, when clarity requires) will denote the Lebesgue measure on $\oi^n$.
We simplify notation by denoting a Borel measure and the state it induces by the same Greek letter.
As remarked in~\cite{mundici95}, $\lambda$ has the following properties:
\begin{itemize}
\item[(A)] $\lambda(f)\in\Qbb$ for every $f\in\Free_n$;
\item[(B)] it is \newword{faithful}: $f\not=0$ implies $\lambda(f)\not=0$;
\item[(C)] it is \newword{automorphism-invariant}: for every $f$ and every automorphism $\sigma$ of $\Free_n$, we have $\lambda(\sigma(f))=\lambda(f)$ (equivalently, for every Borel $A\subseteq\oi^n$ and every McNaughton homeomorphism $S$, we have $\lambda(A)=\lambda(S^{-1}A)$).
\end{itemize}

These properties do not characterize $\lambda$. Indeed, define the
\newword{denominator} of the rational point $p=(\vect r1n)\in\oi^n\cap\Qbb^n$ as the unique integer $d=\den(p)\ge1$ such that $\vect {dr}1n,d$ are relatively prime integers. For a fixed $d$, the set $A^n_d$ of all points in $\oi^n$ having denominator $d$ has finite cardinality $\#A^n_d$, and every McNaughton homeomorphism permutes $A^n_d$.

\begin{example}
Choose\label{ref2} $d\ge1$, and let $\nu^n_d$ be the counting measure on $A^n_d$, giving each point mass $1/\#A^n_d$. Define
$$
\mu^n=\frac{1}{2}\lambda^n+\frac{1}{2}\nu^n_d.
$$
Then $\mu^n$ satisfies 
(A), (B), (C). These measures may exhibit pathological behaviour. Take, e.g., $d=4$ and $f(x_1)=\neg(\neg x_1\oplus\neg x_1)$. We have $A^1_4=\{1/4,3/4\}$, and direct computation shows that integration w.r.t.~$\mu^1$ assigns $f$ value $1/4$. But we may as well consider $f$ as an element of $\Free_2$, possibly making explicit the dummy variable as, say, $f(x_1,x_2)=\neg(\neg x_1\oplus\neg x_1)\land(x_2\oplus\neg x_2)$. Then $\#A^2_4=16$, and $\mu^2(f)=17/64$. The trouble here is that $A^2_4$ does not map onto $A^1_4$ under the canonical projection.
\end{example}

Beyond showing that (A), (B), (C) do not characterize $\lambda$,
Example~\ref{ref2} makes clear that we cannot assign measures to the various $n$-cubes in an unrelated way, and explains the need for the following condition:
\begin{itemize}
\item[(D)] for each $n\ge1$, let $\mu^n$ be a Borel probability measure on $\oi^n$. Let $\pi:\oi^{n+1}\to\oi^n$ be the canonical projection on the first $n$ coordinates. We say that the system $\{\mu^n\}$ is \newword{coherent} if $\mu^n$ and the push-forward measure $\pi_*\mu^{n+1}$ coincide for every~$n$.
\end{itemize}

We now state a condition that will characterize the system of Lebesgue measures.

\begin{definition}
As usual in lattice theory, the \newword{pseudocomplement} of $f\in\Free_n$ is the largest $g$ (if it exists) such that $f\land g=0$. 
We say that $f$ is \newword{pseudotrue} if it has pseudocomplement $0$. Let $\mu$ be a Borel measure on $\oi^n$, and write $kf$ for $f\oplus f\oplus\cdots\oplus f$ ($k$ times). We say that $\mu$ \newword{detects pseudotruths} if the following holds:
\begin{itemize}
\item[(E)] for every pseudotrue $f$, $\sup\{\mu(kf):k\ge1\}=1$.
\end{itemize}
\end{definition}

The \newword{$0$-set} of $f\in\Free_n$ is $Zf=\{p\in\oi^n:f(p)=0\}$. We have that $Zf=\abs{W}$ for some rational polyhedral complex $W$, and conversely any such complex is the $0$-set of some $f$; this follows from the theory of Schauder hats~\cite{mundicimn}, \cite{pantilu}. The sequence $f,2f,3f,\ldots$ is nondecreasing, and its pointwise limit is the characteristic function of the \newword{support} of $f$, the latter being $\supp f=\oi^n\setminus Zf$. For every Borel measure $\mu$ on $\oi^n$, we have $\sup \mu(kf)=\lim \mu(kf)=\mu(\supp f)$ by the Monotone Convergence Theorem.
For every nonempty open set $U$ there exists a $g\not=0$ such that $\supp g\subseteq U$~\cite[Proposition~4.17]{mundicijfa}. Therefore, $f$ is pseudotrue iff $f\land g\not=0$ for every $g\not=0$ iff $\supp g\not\subseteq Zf$ for every $g\not=0$ iff
$Zf$ has empty topological interior.
Summing up, (E) is equivalent to
\begin{itemize}
\item[(E${}'$)] if $W$ is a rational polyhedral complex such that $\abs{W}$ has empty interior, then $\mu(\abs{W})=0$.
\end{itemize}

The following is our main result.

\begin{theorem}
For\label{ref4} every $n\ge1$, let $\mu^n$ be an automorphism-invariant Borel measure on $\oi^n$ that detects pseudotruths. Then, for every $n\ge2$, we have $\mu^n=\lambda^n$. If moreover the system $\{\mu^n\}$ is coherent, then $\mu^1=\lambda^1$ as well.
\end{theorem}

Let us look at Theorem~\ref{ref4} from the viewpoint of \llgroups. Let $P=\{p\in\Rbb^n:\vect p1n\ge0\}$ be the positive cone of $\Rbb^n$, and let $G_n$ be the \llgroup\ of all homogeneous pwl functions with integer coefficients from $P$ to $\Rbb$, with pointwise operations. It is well known that $G_n$ is a projective \llgroup, and all finitely generated projective \llgroups\ have such a form, i.e., are groups of homogeneous pwl functions with integer coefficients defined over some rational cone complex. The map $p\mapsto\{f\in G_n:f(p)=0\}$ is a homeomorphism between the $(n-1)$-dimensional simplex $\Delta^{n-1}=\{\sum \alpha_i e_i:\alpha_i\ge0\text{ and }\sum\alpha_i=1\}$ (the $e_i$'s being the standard basis vectors of $\Rbb^n$) and the maximal spectrum of $G_n$. Now, $G_n$ has plenty of automorphisms, and there is no reason for the existence of a nontrivial probability measure on $\Delta^{n-1}$ invariant under [the duals of] all such automorphisms. The following simple example may clarify this situation.

\begin{example}
Let $n=2$, and let
\begin{align*}
p_1&=e_1+e_2, & p_2&=e_1+2e_2, \\
p'_1&=2e_1+e_2, & p'_2&=e_1+e_2.
\end{align*}
Let $s:P\to P$ be the homogeneous pwl map that fixes $e_1$ and $e_2$, maps $p_i$ to $p'_i$, and is linear on each of the three cones $\Rbb^+e_1+\Rbb^+p_1$, $\Rbb^+p_1+\Rbb^+p_2$, $\Rbb^+p_2+\Rbb^+e_2$.
Then $s$ is a homeomorphism of $P$, and the linear maps on each of the three cones above have integer coefficients. Therefore $s$ is
induced by a unique automorphism $\sigma$ of $G_2$, which is explicitly given by
\begin{align*}
\sigma(x_1)&=x_1\lor\bigl[(3x_1-x_2)\land(x_1+x_2)\bigr], &
\sigma(x_2)&=(-x_1+x_2)\lor\bigl[x_1\land x_2\bigr].
\end{align*}
Parametrizing $\Delta^1$ via $\oi\ni t\mapsto (1-t)e_1+te_2$, one sees that the dual $S$ of $\sigma$ is piecewise-fractional, namely
$$
S(t)=\begin{cases}
t/(1+t), & \text{if $0\le t<1/2$;} \\
(1-t)/(4-5t), & \text{if $1/2\le t<2/3$;} \\
(2t-1)/t, & \text{if $2/3\le t\le1$.}
\end{cases}
$$
\begin{figure}[H]
\begin{center}
\includegraphics[width=4cm,height=4cm]{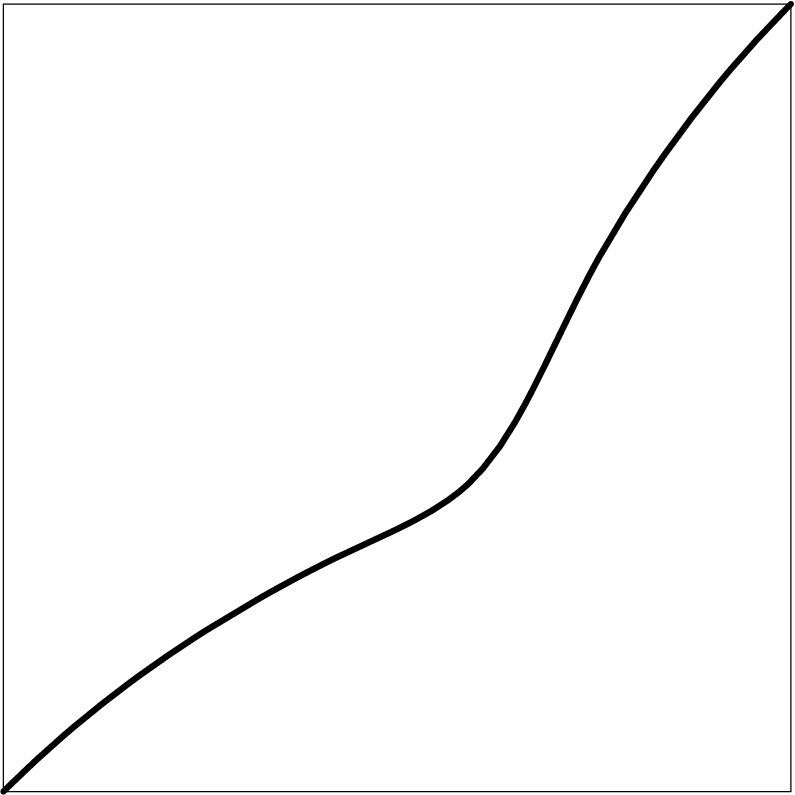}
\end{center}
\end{figure}
All $t\in\oi$, except the repelling fixed point $1$, is attracted under $S$ to the point~$0$. This implies that the only $S$-invariant probability measures on $\oi\simeq\Delta^1$ are the convex combinations of the Dirac measures at $0$ and $1$, a definitely trivial situation.
\end{example}

Things become more interesting if we fix a strong unit $u$ in $G_n$, and restrict attention to those automorphisms of $G_n$ that fix $u$. This amounts to studying the MV-algebra $\Gamma(G_n,u)$ and its group of automorphisms;  the intuition behind is that fixing a strong unit in an \llgroup\ amounts to fixing a ``scale of magnitude''. This intuition is now supported by our Theorem~\ref{ref4}. Indeed, a clever choice of $u$ makes $\Gamma(G_n,u)$ isomorphic to $\Free_{n-1}$: this fact is proved in~\cite[Theorem~4.1]{pantigeneric}, and we leave as an exercise for the reader to show that we can take $u=(x_1\lor x_2\lor\cdots\lor x_{n-1})+x_n$. Then the Lebesgue measure $\lambda$ on the 
maximal spectrum $\Delta^{n-1}$ of $G_n$ is invariant under 
[the duals of] all automorphisms of $G_n$ that leave $u$ fixed, and
Theorem~\ref{ref4} says that this property, together with being null 
on underdimensioned $0$-sets, does indeed characterize $\lambda$.

Finally, the relevance of Theorem~\ref{ref4} from the viewpoint of algebraic logic is clear: it gives a distinguished status to the Lebesgue measure as an averaging measure on the space $\oi^n$ of truth-value assignments in \Luk\ infinite-valued logic. We only stress here that, although the proof of Theorem~\ref{ref4} exploits a fair amount of the arithmetic and piecewise-linear structure of $\oi^n$, the resulting characterization of the state on $\Free_n$ induced by integration w.r.t.~the Lebesgue measure as the only automorphism-invariant state detecting pseudotruths is intrinsic and purely algebraic. We refer to~\cite[\S2]{pantibernoulli} for a detailed discussion on this issue.

\section{The case $n=2$}

The key point in establishing Theorem~\ref{ref4} is the case $n=2$.

\begin{theorem}
Let\label{ref6} $\mu$ be a Borel measure on $\oi^2$. Assume that $\mu$
detects pseudotruths and is automorphism-invariant. Then $\mu=\lambda$.
\end{theorem}

The analogous statement is false for $n=1$. Indeed, it is easy to show that the only McNaughton homeomorphisms of $\oi$ are the identity and the flip $p\mapsto 1-p$. Therefore, any flip-invariant strictly positive density $s$ (i.e., a $\lambda$-integrable function $s:\oi\to(0,\infty)$ such that $\int s\,d\lambda=1$ and $s(p)=s(1-p)$) determines a measure absolutely continuous w.r.t.~$\lambda$, automorphism-invariant, faithful, and detecting pseudotruths. The trouble here is that $\oi$ has so few McNaughton homeomorphisms that rigidity phenomena cannot occur.

Rigidity is a recurrent theme in measurable dynamics. If $X$ is a space on which a transformation $R$ acts, then typically there is a multitude of $R$-invariant probability measures on $X$. But if we add a second transformation $S$ which is incommensurable with $R$ (i.e., no power of $R$ is a power of $S$), then the set of measures invariant for both $R$ and $S$ may shrink down drastically, often just to the ``obvious'' ones. The most famous, and still unsolved, problem in this area is the Furstenberg Conjecture: the only nonatomic probability on the unit circle $\{z\in\Cbb:\abs{z}=1\}$ which is invariant under both $z\mapsto z^2$ and $z\mapsto z^3$ is the Lebesgue measure.

We shall prove Theorem~\ref{ref6} using such a rigidity approach.
I must sincerely thank Carlangelo Liverani and Fran\c{c}ois Ledrappier for explaining me the rigidity of the two-dimensional torus under the action of $\SL_2\Zbb$; the insight they gave was the starting point of the present work.

We will denote rational points of the unit square switching freely between nonhomogeneous and homogeneous coordinates, e.g., the center of the square is either $(1/2,1/2)$ or $(1:1:2)$. For $h\ge0$, define 
$p^h_0=(h+1:h:2h+1)$, $p^h_1=(h+1:h+1:2h+1)$, $p^h_2=(h:h+1:2h+1)$, $p^h_3=(h:h:2h+1)$. For every $k\ge0$, we define a McNaughton homeomorphism $R_k$ by considering the following rational polyhedral complex (we draw the picture for $k=1$):

\begin{figure}[!h]
\begin{center}
\includegraphics[width=4cm]{figura-1}
\end{center}
\end{figure}

The vertices of the intermediate square are $p^k_i$, for $i=0,\ldots,3$, and those of the inner square are $p^{k+1}_i$. Let $R'_k$ be the unique homeomorphism that:
\begin{enumerate}
\item maps $p^{k+1}_i$ to $p^{k+1}_{i+1\pmod{4}}$;
\item fixes every other vertex;
\item is affine linear on each polyhedron.
\end{enumerate}
By direct computation, one sees that $R'_k$ is a McNaughton homeomorphism. E.g., if $T$ is the triangle $\angles{p^{k+1}_0,p^k_0,p^k_1}$, then $R'_k[T]$ is the triangle
$\angles{p^{k+1}_1,p^k_0,p^k_1}$, and the matrix expressing $R'_k\restriction T$ in homogeneous coordinates is
\begin{equation}\tag{$*$}
\begin{pmatrix}
1 & 0 & 0 \\
-2k-1 & 1 & k+1 \\
0 & 0 & 1
\end{pmatrix}
\end{equation}
Define $R_k=(R'_k)^4$ to be the fourth power of $R'_k$.

Let now $D$ be the unit disk in polar coordinates $(r,\theta)$, for $-\pi<\theta\le\pi$, and let $M:D\to\oi^2$ be the homeomorphism defined by:
$$
M(r,\theta)=(1/2,1/2)+\begin{cases}
r(-1/2, -2\theta/\pi-3/2), & \text{if $-\pi<\theta\le-\pi/2$;} \\
r(2\theta/\pi+1/2,-1/2), & \text{if $-\pi/2<\theta\le0$;} \\
r(1/2,2\theta/\pi-1/2), & \text{if $0<\theta\le\pi/2$;} \\
r(-2\theta/\pi+3/2,1/2), & \text{if $\pi/2<\theta\le\pi$.}
\end{cases}
$$
For every $r\in\oi$, let $K_r=\{(r,\theta):-\pi<\theta\le\pi\}$ be the circle of radius $r$, let $Q_r=M[K_r]$ (it is the boundary of a square centered in $(1/2,1/2)$), and let $t_k:\oi\to\oi$ be the continuous piecewise-fractional function defined by
$$
t_k(r)=\min(1,\max\left(0,1/(2r)-k-1/2)\right).
$$

\begin{lemma}
The\label{ref7} homeomorphism $M$ provides a conjugation between $R_k$ and the map $T_k:D\to D$ given by $T_k(r,\theta)=(r,\theta+2\pi t_k(r))$ (i.e., $R_k\circ M=M\circ T_k$).
\end{lemma}
\begin{proof}
Let $v\in K_r$. We may assume $1/\bigl(2(k+1)+1\bigr)<
r<1/\bigl(2k+1\bigr)$, since otherwise $T_k$ fixes $v$ and $R_k$ fixes $M(v)$. If $u$ is another point in $K_r$, we define the \newword{oriented distance} between $u$ and $v$ as the length $d(u,v)$ of the path which goes from $u$ to $v$ along $K_r$ in the counterclockwise direction (we normalize the total length of $K_r$ to be $1$). We define analogously the oriented distance of two points of $Q_r$. The differential of $M$ maps the unit tangent vector $\partial/\partial\theta\rvert_u$ to $K_r$ at $u$ to a vector of constant modulus, namely either
$-2r/\pi\,\partial/\partial y\rvert_{Mu}$, or
$2r/\pi\,\partial/\partial x\rvert_{Mu}$, or
$2r/\pi\,\partial/\partial y\rvert_{Mu}$, or
$-2r/\pi\,\partial/\partial x\rvert_{Mu}$, depending on $\theta\in(-\pi,\pi)\setminus\{0,\pm \pi/2\}$.
Therefore, if $u,v,u',v'\in K_r$ are such that $d(u,v)=d(u',v')$, then $d(Mu,Mv)=d(Mu',Mv')$. Both $T_k$ and $R_k$ preserve the oriented distance: this is clear for $T_k$ and can easily be checked for $R'_k$, whence for $R_k$. Fix $u=(r,0)$; explicit computation shows that
$$
u\overset{M}{\longmapsto}(1+r:1-r:2)
\overset{R'_k}{\longmapsto}(1+r:-2kr-2r+2:2)
\overset{M^{-1}}{\longmapsto}(r,\pi/2\,t_k(r)),
$$
(the second step in the computation uses the matrix ($*$) displayed above), and therefore $R_kMu=(R'_k)^4Mu=MT_ku$. 
We thus obtain $d(MT_ku,R_kMv)=d(R_kMu,R_kMv)=d(Mu,Mv)=d(MT_ku,MT_kv)$ (because $d(u,v)=d(T_ku,T_kv)$). Therefore $R_kMv$ and $MT_kv$ have the same distance from $MT_ku$, whence they are equal.
\end{proof}

Until the end of this section $\mu$ is a fixed measure satisfying the hypotheses of Theorem~\ref{ref6}. The function $\rho:\oi^2\to\oi$ that maps all points in $Q_r$ to $r$ is Borel measurable. Hence the family $\xi=\{Q_r:r\in\oi\}$ of its fibers is a measurable decomposition of $\oi^2$, and $\mu$ disintegrates over $\xi$~\cite[p.~26]{rohlin49a}. This means the following: write $\nu=\rho_*\mu$ for the push-forward of $\mu$ on $\oi$ via $\rho$. Then for every $r$ there exists a Borel probability measure $\mu_r$ on $Q_r$ such that, for every $\mu$-measurable $A$, we have:
\begin{itemize}
\item[(a)] $A\cap Q_r$ is $\mu_r$-measurable;
\item[(b)] the map $r\mapsto\mu_r(A\cap Q_r)$ is $\nu$-measurable;
\item[(c)]
$$
\mu(A)=\int_{\oi}\mu_r(A\cap Q_r)\,d\nu;
$$
\item[(d)] if $\{\mu'_r\}$ is another system of conditional measures satisfying the above properties, then $\mu'_r=\mu_r$ for $\nu$-every $r$.
\end{itemize}

\begin{lemma}
For\label{ref11} $\nu$-every $r$, the conditional measure $\mu_r$ coincides with the linear Lebesgue probability measure $\tau_r$ on $Q_r$.
\end{lemma}
\begin{proof}
It is immediate that $Q_r$ contains a rational point iff the intersection point $(1+r:1-r:2)$ of $Q_r$ with the line segment $\angles{(1:1:2),p^0_0}$ is a rational point iff $r$ is a rational number.
Since $\mu$ satisfies (E${}'$), we have
$$
\nu(\oi\cap\Qbb)=\mu\bigl(\bigcup\{Q_r:r\in\oi\cap\Qbb\}\bigr)=0,
$$
whence we can neglect rational $r$'s. Fix $k\ge0$; it suffices to show that the set of all irrational $r$ in the interval $(1/(2(k+1)+1),1/(2k+1))$ such that $\mu_r\not=\tau_r$ has $\nu$-measure $0$. Given such an $r$, consider the topological dynamical system $(K_r,T_k\restriction K_r)$: it is a rotation by an angle of $(1/r-2k-1)\pi$, which is an irrational multiple of $\pi$. 
Irrational rotations are uniquely ergodic, i.e., preserve a unique probability Borel measure~\cite[\S~6.5]{Walters82}. By Lemma~\ref{ref7}, $(K_r,T_k\restriction K_r)$ is topologically conjugate to $(Q_r,R_k\restriction Q_r)$, which is therefore 
uniquely ergodic as well. Since $R_k\restriction Q_r$ preserves the oriented distance, it preserves $\tau_r$. On the other hand, by~\cite[\S~2.5]{rohlin49b}, $R_k\restriction Q_r$ preserves $\mu_r$ for $\nu$-every $r$; this settles our claim.
\end{proof}

For $i=0,\ldots,3$, let $T_i$ be the open triangle
$$
T_i=M\bigl[\{(r,\theta):0<r<1, i\pi/2<\theta<(i+1)_{(\mmod 4)}\pi/2\}\bigr].
$$

\begin{lemma}
Let\label{ref8} $A$ be a Borel subset of the unit square, let $e_1=(1,0)$, $e_2=(0,1)$ be the standard basis of $\Rbb^2$, and let $s\in\Rbb$. If $A,se_1+A\subseteq T_1\cup T_3$ then $\mu(A)=\mu(se_1+A)$ and, analogously, if 
$A,se_2+A\subseteq T_0\cup T_2$ then $\mu(A)=\mu(se_2+A)$.
\end{lemma}
\begin{proof}
Since the measures $\mu_r$ are linear Lebesgue measures, they are invariant under ``horizontal'' translations inside $T_1\cup T_3$, and under ``vertical'' ones inside $T_0\cup T_2$. The statement is now immediate from the disintegration property~(c).
\end{proof}

\begin{lemma}
Let\label{ref5} $T\subset\Rbb^n$ be a compact convex polyhedron of dimension $n$, $T^o$ its topological interior, $\eta$ a finite Borel measure on $T^o$. Let us assume that for every $v\in\Rbb^n$ and every Borel set $A$ such that $A,v+A\subseteq T^o$, we have $\eta(A)=\eta(v+A)$. Then there exists $c>0$ such that $\eta=c\lambda$.
\end{lemma}
\begin{proof}
The proof is modeled on that of~\cite[Theorem~2.20(d)]{rudin87}, so we just sketch it. For $k\ge0$, a \newword{$k$-box} is a translate of $[0,2^{-k})^n$ by a point in $2^{-k}\Zbb^n$. Since the boxes contained in $T^o$ generate the Borel $\sigma$-algebra and $\eta,\lambda$ are regular, it suffices to find $c$ such that $\eta(B)=c\lambda(B)$ for every box $B\subset T^o$. Let $h$ be minimum such that there exists an $h$-box $C$ in $T^o$, and let $c=\eta(C)/\lambda(C)$. Let now $k\ge h$, and let $B$ be any $k$-box in $T^o$. Then $C$ is partitioned in $2^{n(k-h)}$ $k$-boxes and each of them, by translation invariance, has the same $\eta$-measure and the same $\lambda$-measure as $B$. Therefore
$$
\eta(B)=\frac{\eta(C)}{2^{n(k-h)}}=\frac{c\lambda(C)}{2^{n(k-h)}}=
c\lambda(B).
$$
\end{proof}

In order to apply Lemma~\ref{ref5}, we need to find another family $\{S_k:k\ge0\}$ of McNaughton homeomorphisms that force $\mu$ to be invariant under translates along new (i.e., neither horizontal nor vertical) directions. For $k\ge0$, we consider the following complexes (again, we drawn a picture for $k=1$):

\begin{figure}[!h]
\begin{center}
\includegraphics[width=4cm]{figura-2}
\end{center}
\end{figure}

Let $q^h_0=(h+2:h+1:2h+2)$, $q^h_1=(h+1:h+2:2h+2)$,
$q^h_2=(h:h+1:2h+2)$, $q^h_3=(h+1:h+2:2h+2)$. The vertices of the exterior rhombus are then $q^0_i$, for $i=0,\ldots,3$, those of the intermediate one are $q^k_i$, and those of the inner one are $q^{k+1}_i$. Let $S'_k$ be the unique homeomorphism that:
\begin{enumerate}
\item maps $q^{k+1}_i$ to $q^{k+1}_{i+1\pmod{4}}$;
\item fixes every other vertex;
\item is affine linear on each polyhedron.
\end{enumerate}
In this case, $S'_k$ is not a McNaughton homeomorphism, since the matrices that express the map have rational non-integer entries. Nevertheless, exactly the same proof as in~\cite[Lemma~3.1]{pantibernoulli} shows that the square of $S'_k$ is a McNaughton homeomorphism; a fortiori, $S_k=(S'_k)^4$ is a McNaughton homeomorphism.

Define $N:D\to\oi^2$ by
$$
N(r,\theta)=(1/2,1/2)+\begin{cases}
r(\theta/\pi+1/2, -\theta/\pi-1), & \text{if $-\pi<\theta\le-\pi/2$;} \\
r(\theta/\pi+1/2,\theta/\pi), & \text{if $-\pi/2<\theta\le0$;} \\
r(-\theta/\pi+1/2,\theta/\pi), & \text{if $0<\theta\le\pi/2$;} \\
r(-\theta/\pi+1/2,-\theta/\pi+1), & \text{if $\pi/2<\theta\le\pi$.}
\end{cases}
$$
Then $N$ is a homeomorphism between $D$ and the rhombus $E=\angles{q^0_0,q^0_1,q^0_2,q^0_3}$. Let $V_r=N[K_r]$, and define $h_k:\oi\to\oi$ by
$$
h_k(r)=\min(1,\max(0,1/(2r)-k-1)).
$$
Explicit computation, exactly as in the proof of Lemma~\ref{ref7}, shows that $N$ provides a conjugation between the restriction of $S_k$ to $E$ and the twist map $(r,\theta)\mapsto(r,\theta+2\pi h_k(r))$ on $D$. Again, $\mu$ disintegrates along the $V_r$'s, and all the argument leading to Lemma~\ref{ref8} carries through.
Setting, for $i=0,\ldots,3$,
$$
H_i=N\bigl[\{(r,\theta):0<r<1, i\pi/2<\theta<(i+1)_{(\mmod 4)}\pi/2\}\bigr],
$$
we obtain the following Lemma.

\begin{lemma}
Let\label{ref9} $A,e_1,e_2,s$ be as in the statement of Lemma~\ref{ref8}. Then:
\begin{itemize}
\item[(i)] if $A,s(e_1+e_2)+A\subseteq H_1\cup H_3$ then $\mu(A)=\mu(s(e_1+e_2)+A)$;
\item[(ii)] if $A,s(-e_1+e_2)+A\subseteq H_0\cup H_2$ then $\mu(A)=\mu(s(-e_1+e_2)+A)$;
\item[(iii)] on $E$, the measure $\mu$ is a positive multiple $c\lambda$ of the Lebesgue measure.
\end{itemize}
\end{lemma}
\begin{proof}
(i) and (ii) are proved exactly as in Lemma~\ref{ref8}.
By Lemma~\ref{ref5}, Lemma~\ref{ref8} and (i), (ii) above, on each open triangle $T_i\cap H_j$ the measure $\mu$ is a positive multiple $c_{ij}\lambda$ of Lebesgue. Every $T_i\cap H_j$ is mapped to any other
$T_{i'}\cap H_{j'}$ by an appropriate symmetry of the square, and these symmetries are McNaughton homeomorphisms. Hence $c_{ij}=c_{i'j'}$  and~(iii) follows since, by Condition~(E${}'$),
$\mu(E\setminus\bigcup\{T_i\cap H_j:i,j=0,\ldots,3\})=0$.
\end{proof}

Each $T_i\setminus E$ is the disjoint union of two open triangles, and each of these triangles is the union (modulo $\mu$-null sets, again by~(E${}'$)) of countably many open rhombi, as in the picture below:

\begin{figure}[!h]
\begin{center}
\includegraphics[width=6cm]{figura-3}
\end{center}
\end{figure}

Let $F$ be one of these rhombi, say for definiteness' sake inside the triangle $\angles{p^0_3,q^0_3,(1/4,1/4)}$. Let $t$ be the horizontal translation that maps $F$ inside $E$ (the grey region above). Then $\mu$ is invariant under translations inside $F$, because any such translation is conjugate via $t$ to a translation inside $E$. By Lemma~\ref{ref5}, $\mu=c_F\lambda$ on $F$, and necessarily $c_F$ equals the constant $c$ in Lemma~\ref{ref9}(iii), because of Lemma~\ref{ref8} and the fact that $F,t[F]\subseteq T_3$.
The union of $E$ with all these countably many rhombi is ---modulo $\mu$-null sets--- the whole of $\oi^2$, and therefore $c=1$. This concludes the proof of Theorem~\ref{ref6}.

\section{Proof of Theorem~\ref{ref4}}

Let $\mu^n$ be measures as in the statement of Theorem~\ref{ref4}. If the system $\{\mu^n\}$ is coherent, then $\mu^n$ determines $\mu^m$ for $1\le m\le n$. By Theorem~\ref{ref6} $\mu^2=\lambda^2$, and hence $\mu^1=\lambda^1$; this settles the last claim in Theorem~\ref{ref4}. If the system is not coherent, then $\mu^1$ may be different from $\lambda^1$, as discussed after the statement of Theorem~\ref{ref6}.

Fix now $n>2$. The natural projection $\pi:\oi^n\ni(\vect x1n)\mapsto
(\vect x3n)\in\oi^{n-2}$ is continuous, whence the set $\xi$ of its fibers is a measurable partition of the $n$-cube. As described after Lemma~\ref{ref7}, $\mu^n$ disintegrates over $\xi$. Let $\nu=\pi_*\mu^n$ be the push-forward of $\mu^n$ on $\oi^{n-2}$, and let $\mu^n_u$ be the conditional measure on $C_u=\pi^{-1}\{u\}\in\xi$.
The group $G_2$ of all McNaughton homeomorphisms of $\oi^2$ embeds in the group $G_n$ of all McNaughton homeomorphisms of $\oi^n$ in the obvious way: if $R\in G_2$, $p\in\oi^2$, and $u\in\oi^{n-2}$, then we set $R(p,u)=(R(p),u)$.

\begin{lemma}
For\label{ref10} every $R\in G_2$, every rational polyhedral complex $W$ in $\oi^2$ such that $\abs{W}$ has empty interior, and
$\nu$-every $u\in\oi^{n-2}$, the following hold:
\begin{itemize}
\item[(i)] $R\restriction C_u$ preserves $\mu^n_u$;
\item[(ii)] $\mu^n_u((\abs{W}\times\oi^{n-2})\cap C_u)=0$.
\end{itemize}
\end{lemma}
\begin{proof}
By construction, every $C_u$ is invariant under the action of $G_2$. For a fixed $R\in G_2$, the set $U_R$ of all $u\in\oi^{n-2}$ such that $\mu_u$ is invariant under $R\restriction C_u$ has $\nu$ measure $1$. Let $W$ be a complex as stated. Then $\abs{W}\times\oi^{n-2}=\abs{V}$ for some rational polyhedral complex $V$ on $\oi^n$. By condition~(E${}'$), $\mu^n(\abs{V})=0$, and therefore the characteristic function of $\abs{V}$ equals $0$ as an element of $L_1(\oi^n,\mu^n)$. By the construction of the conditional measures (see, e.g.,~\cite{rao93}), the set $U'_W$ of all $u$ such that $\mu^n_u(\abs{V}\cap C_u)=0$ has $\nu$ measure $1$. Since $G_2$ is a countable group and there exist countably many complexes $W$ as above, the set of all $u$ that satisfy both~(i) and~(ii) has $\nu$ measure~$1$.
\end{proof}

By Lemma~\ref{ref10} and Theorem~\ref{ref6}, for $\nu$-every $u$ the conditional measure $\mu^n_u$ is the Lebesgue measure 
$\lambda^2$, which is invariant under translations inside $C_u$. The disintegration property
$$
\mu^n(A)=\int_{\oi^{n-2}} \mu^n_u(A\cap C_u)\,d\nu
$$
assures then that $\mu^n$ is invariant under translations of the form $(p,u)\mapsto(q+p,u)$, for $q\in\Rbb^2$. But of course all the above construction can be repeated for any coordinate pair in $\oi^n$, and therefore $\mu^n$ is invariant under any translation inside~$\oi^n$. By Lemma~\ref{ref5}, $\mu^n$ is the Lebesgue measure $\lambda^n$.

\end{document}